\documentclass[11pt]{elsarticle}
\usepackage{lineno}
\modulolinenumbers[1]
\usepackage{hyperref} 

\usepackage[letterpaper, margin=1.0in]{geometry}
\renewcommand{\baselinestretch}{1.2}

\usepackage{soul}

\usepackage{lipsum}
\usepackage{rotating}
\usepackage{pstricks, pst-node, psfrag}
\usepackage{amssymb,amsmath}
\usepackage{mathtools}
\usepackage{verbatim,enumerate}
\usepackage{rotating, lscape}
\usepackage{setspace}
\usepackage[hang, flushmargin]{footmisc}
\usepackage{caption}
\usepackage{subcaption}
\usepackage{cancel}
\usepackage{float}
\usepackage{ragged2e}
\usepackage{centernot}
\usepackage{tikz}
\usepackage{tikzsymbols}
\usetikzlibrary{shapes,arrows,positioning,plotmarks}
\usepackage{textcomp} 
\usepackage{gensymb} 

\tikzstyle{endpt} = [rectangle, draw, fill=red!20,
    text width=9.8em, text centered, rounded corners, minimum height=4em]
\tikzstyle{block} = [rectangle, draw, top color=white, bottom color=blue!20,
    text width=9.8em, text centered, rounded corners, minimum height=4em]
\tikzstyle{line} = [draw, -latex', very thick]

\usepackage{listings}
\usepackage{color} 
\definecolor{mygreen}{rgb}{28,172,0} 
\definecolor{mylilas}{rgb}{170,55,241}
\definecolor{mygray}{rgb}{0.5,0.5,0.5}
\definecolor{mycyan}{rgb}{0,255,255}
\definecolor{magenta}{rgb}{1,0,1}

\definecolor{backgreen}{rgb}{0.00, 0.169, 0.212}
\definecolor{textgray}{rgb}{0.514, 0.580, 0.589}

\usepackage[T1]{fontenc}
\newcommand*\patchAmsMathEnvironmentForLineno[1]{%
    \expandafter\let\csname old#1\expandafter\endcsname\csname #1\endcsname
    \expandafter\let\csname oldend#1\expandafter\endcsname\csname end#1\endcsname
    \renewenvironment{#1}%
        {\linenomath\csname old#1\endcsname}%
        {\csname oldend#1\endcsname\endlinenomath}}%
\newcommand*\patchBothAmsMathEnvironmentsForLineno[1]{%
    \patchAmsMathEnvironmentForLineno{#1}%
    \patchAmsMathEnvironmentForLineno{#1*}}%
\AtBeginDocument{%
    \patchBothAmsMathEnvironmentsForLineno{equation}%
    \patchBothAmsMathEnvironmentsForLineno{align}%
    \patchBothAmsMathEnvironmentsForLineno{flalign}%
    \patchBothAmsMathEnvironmentsForLineno{alignat}%
    \patchBothAmsMathEnvironmentsForLineno{gather}%
    \patchBothAmsMathEnvironmentsForLineno{multline}%
}

\newcommand{\bl}{\begin{linenomath}}
\newcommand{\el}{\end{linenomath}}
\newcommand{\be}{\begin{equation}}
\newcommand{\ee}{\end{equation}}

\newcommand{\dt}{\Delta t}

\renewcommand{\hat}{\widehat}

\DeclareMathAlphabet{\mathpzc}{OT1}{pzc}{m}{it}
\DeclareMathOperator\SSE{\mathrm{SSE}}
\DeclareMathOperator\mcE{\mathcal{E}}
\DeclareMathOperator*\argmin{\mathrm{arg \ min}}

\newcommand{\rev}[1]{#1}

\journal{Advances in Water Resources}

\biboptions{sort&compress}
\begin{document}

\begin{frontmatter}
\title{A Computational Information Criterion for Particle-Tracking with Sparse or Noisy Data\tnoteref{mytitlenote}}
\tnotetext[mytitlenote]{This work was supported by the National Science Foundation under awards DMS-1614586 and DMS-1911145 and the US Army Research Office under Contract/Grant number W911NF-18-1-0338.}

\author{Nhat Thanh Tran\fnref{uci}}
\ead{nhattt@uci.edu}
\author{David A. Benson\fnref{hydro}}
\ead{dbenson@mines.edu}
\author{Michael J. Schmidt\fnref{snl}}
\ead{mjschm@sandia.gov}
\author{Stephen D. Pankavich\fnref{ams}}
\ead{pankavic@mines.edu}

\fntext[uci]{Department of Mathematics, University of California, Irvine, Irvine, CA, 92697, USA}
\fntext[hydro]{Hydrologic Science and Engineering Program, Department of Geology and Geological Engineering, Colorado School of Mines, Golden, CO, 80401, USA}
\fntext[snl]{Center for Computing Research, Sandia National Laboratories, Albuquerque, NM 87185, USA}
\fntext[ams]{Department of Applied Mathematics and Statistics, Colorado School of Mines, Golden, CO, 80401, USA}

\begin{abstract}
Traditional probabilistic methods for the simulation of advection-diffusion equations (ADEs) often overlook the entropic contribution of the discretization, e.g., the number of particles, within associated numerical methods. Many times, the gain in accuracy of a highly discretized numerical model is outweighed by its associated computational costs or the noise within the data. We address the question of how many particles are needed in a simulation to best approximate and estimate parameters in one-dimensional advective-diffusive transport. To do so, we use the well-known Akaike Information Criteri\rev{on} (AIC) and a recently-developed correction called the Computational Information Criteri\rev{on} (COMIC) to guide the model selection process. Random-walk and mass-transfer particle tracking methods are employed to solve the model equations at various levels of discretization. Numerical results demonstrate that the COMIC provides an optimal number of particles that can describe a more efficient model in terms of parameter estimation and model prediction compared to the model selected by the AIC even when the data is sparse or noisy, the sampling volume is not uniform throughout the physical domain, or the error distribution of the data is non-IID Gaussian.
\end{abstract}


\end{frontmatter}

\newpage
\section{Introduction}
Numerical methods of all sorts are used to approximate the solutions of various model equations in hydrology.  The independent variables in these models are discretized, and model coefficients are populated so as to faithfully reproduce some set of measured dependent variables (i.e., data).  Of course, both the model solution and the measured data will contain errors; therefore, a perfect match of model to data is not necessarily desired.  In fact, a perfect match will most often reduce the ability of the model to predict new sets of data, because the model is overfit to a single realization of noise (see the excellent overview by \cite{Konishi_book}). \rev{The bias between a model that is overfit and an underlying ``true'' model was classically addressed by Akaike \cite{Akaike1974} by considering approximate measures of the entropy of the probability distributions associated with the likelihood that data arises from a specific candidate model}.  In short, these measures introduce an entropic penalty to maximum likelihood estimates of goodness of fit when the number of {\em adjustable} parameters or coefficients increases.  In a similar fashion, the computational entropy of a model is a bias against its fitness and needs to be accounted for when using a model in a predictive mode.

For a few examples, consider first a particle-tracking model of contaminant transport.  It is often visually pleasing to use a large number of particles in order to obtain a smooth (i.e., low \rev{concentration} variance) histogram for comparison with data.  However, if the noise in the data far exceeds the noise in the histogram, then the large number of particles is superfluous from a model prediction perspective. Indeed, the model smoothness is not implied by the data and should not be included in a prediction.  In addition, if the model is used repetitively in a parameter estimation procedure, then these extra calculations may also become tangibly burdensome. Another example is stochastic Monte Carlo modeling of groundwater flow using an Eulerian (e.g., finite-difference) simulator.  Oftentimes the spatial and temporal discretizations are thought to be a free modeling choice, and random realizations of hydraulic conductivity are generated, yielding models with arbitrary degrees of freedom.  Subsequently, the models are weighted in relation to the goodness of fit (e.g., \cite{Bevin_Binley_glue,Poeter_2005,Ye_2004}), but one must ask whether simpler models, in terms of {\em both} discretization and parameterization, should be elevated because of their computational simplicity.  

Recently, \cite{Entropy} showed that a computational entropy penalty\rev{, called the Computational Information Criterion or COMIC, can be easily used to address this issue for} simple systems (i.e., those with constant discretization or particle numbers), but it has yet to be shown that those results can be generalized to more realistic modeling scenarios.  Here\rev{,} we show that measures of computational (entropy) penalty may be extended to less ideal cases.  
\rev{In the next section, we will review both the COMIC and the advection-dispersion (mixed hyperbolic and parabolic) partial differential equation (PDE) model of interest, and then discuss associated Lagrangian numerical methods. Sections \ref{sec:sparse} and \ref{sec:nonuniform} are devoted to understanding the effects of sparse and non-uniformly spaced data sets, respectively, on the information criterion and subsequent parameter estimation. Next, we reformulate the COMIC in Section \ref{sec:dv} to account for a non-uniform discretization volume and demonstrate its equivalence (in terms of particle number selection) with the COMIC derived from a uniform sampling volume. In Section \ref{sec:noGauss}, we derive a more general version of the COMIC for non-Gaussian errors and non-uniform error variance and test it on noisy data sets. Conclusions are discussed and summarized in the final section. }

\rev{\section{Model and Methods}}
\label{sec:mm}
To enable a complete analysis, we use an equation that does not require a numerical solution (i.e., an analytical solution is readily available).  This restriction is not required.  The simplest case of the one-dimensional, constant-coefficient advection-diffusion equation (ADE) is given by
\begin{equation}\label{eq:ADE}
\dfrac{\partial c}{\partial t} = -v \dfrac{\partial c}{\partial x} + D \dfrac{\partial^2 c}{\partial x^2} \, ,
\end{equation}
where $c(x,t)$ is the solution to the PDE, $v$ is \rev{the velocity}, and $D$ is the diffusion coefficient. \rev{As we are interested in particle methods}, we will restrict our attention to the initial condition $c(x,0)=\delta(x_0)$ so that $c(x,t)$ is a probability density function (PDF).  In practice, the solution to the PDE is approximated by the numerical method $c_n(x,t)$, which is a function of a discretization parameter $n$, representing the number of particles in a particle method or nodes in a finite-difference approximation. The choice of $n$ is arbitrary and we seek a systematic way to choose this modeling parameter.

The newly developed COMIC may be used to select models amongst these different levels of discretization \cite{Entropy}. This \rev{criterion} is a extension of Akaike's ``an information criterion'' (AIC) in which there is a penalization of the usage of more information (in the form of adjustable parameters) to describe the model.
For completeness, we recall that the AIC is defined by
\begin{equation}\label{eq:AIC}
    \mathrm{AIC} = -2\ln (\mathcal{L}(\hat{\theta}))  + 2p
\end{equation}
where $\mathcal{L}(\hat{\theta})$ is the likelihood function evaluated at the maximum likelihood estimate for the unknown parameters $\theta$, and $p$ is the number of parameters.
For a computational model, the number of nodes or particles contributes to the entropy and must be accounted for when comparing model predictive fitness.  In short, the Kullback-Leibler \cite{K_L_1951,Kullback_book} relative entropy for a discretized model contributes an extra term to the AIC, and the COMIC takes the form
\begin{equation}\label{eq:COMIC_exact}
    \mathrm{COMIC} = \mathrm{AIC}  - \int c(x,t)\ln(\Delta V(x))\; dx
\end{equation}

\noindent where $\Delta V(x)$ is the sampling volume, here given by the spacing between the particles or nodes. In the case that $\Delta V(x) = |\Omega|/ n$ is constant where $|\Omega|$ represents the length of the spatial domain, then the COMIC will merely reduce to
\begin{equation}\label{eq:COMIC}
    \mathrm{COMIC} = \mathrm{AIC}  + \ln(n)
\end{equation}
up to a constant.
It is this quantity that should be minimized in order to identify the best available predictive model \cite{Konishi_book,Entropy}.

For this study, we use two Lagrangian methods: random-walk particle-tracking (RWPT) and mass-transfer particle-tracking (MTPT).  Both methods are described in detail elsewhere \cite{labolle,Salamon_2006,Benson_arbitrary,Schmidt_accuracy}, so we only summarize here.  \rev{In its simplest form}, the RWPT method places a number of particles $n$ at the release position $x_0$, each with constant mass $1/n$.  At every time step of duration $\Delta t$, each particle moves with mean $v\Delta t$ and random deviation $\sqrt{2D\Delta t} \xi$, where $\xi \sim \mathcal{N}(0,1)$ is a standard Normal random variable.  At any desired time, bins of size $\Delta x_i$, centered at points $x_i$ are constructed and the particle count $n_i$ in each bin is converted to concentration by $c_n(x_i,t)= n_i/(n\Delta x_i)$.

\rev{Contrastingly, the MTPT method typically} allows a portion of the dispersion to be performed by random walks as above, and the remaining portion is performed by mass transfer between particles \rev{\cite{Benson_Poise,Entropy,guillem_MRIEM,herrera_2009,Engdahl_WRR,Schmidt_fluid_solid}}.
The mass transfer between any and all particles is governed by the equation
\begin{equation}
\label{eqn:MTPT}
m_i(t+\dt)=m_i(t)-\sum_{i=1}^n \left (m_i(t)-m_j(t) \right) W_{ij},
\end{equation}
where for each particle pair denoted $i,j$,
\begin{equation}
W_{ij}=\frac{(1/\sqrt{4\pi D\Delta t})\exp(-s_{ij}^2/4D\Delta t)}{\rho_{ij}}
\end{equation}
is the normalized kernel that determines the weight of mass transfer between particles $i$ and $j$ \cite{Benson_arbitrary}, $\rho_{ij}$ is a normalizing constant \rev{that ensures conservation of mass and is typically taken to be the particle density \cite{Guillem_SPH,Schmidt_discon}}, and $s_{ij}$ is the distance between particles $i$ and $j$.
We note that the choice of the Gaussian kernel's bandwidth is a free parameter\rev{. In this case we choose it to be $\sqrt{2 D \dt}$}, making $W$ a normalized version of the fundamental solution to the diffusion equation\rev{, and this is equivalent to choosing $\beta = 1$ within the convention of \cite{Guillem_SPH}}.
Because the particles continually change mass, they are typically pre-distributed throughout the domain with zero mass.  One particle at the release point is given unit mass.  For this study, when we employ an MTPT method all dispersion is modeled by mass transfer, i.e., the particles move only by their mean velocity with no random motion (in contrast to the RWPT simulations that only simulate dispersion via random-walk, with no mass transfer). 

In subsequent sections we perform a variety of numerical simulations using the COMIC to select the optimal number of particles in a model and perform parameter estimation.
All numerical simulations were conducted in MATLAB using a desktop computer with a 3.5 GHz Intel Core i7 processor and 16 GB of RAM. The code used to generate the results in this section is available at \url{https://github.com/nhatthanhtran/entropy2020} \cite{EntropyTech_repo}.


\section{Sparsity of Data}
\label{sec:sparse}
\begin{figure}
\centering
\hspace{-0.35in}
\begin{subfigure}[b]{\textwidth}
    \centering
    \includegraphics[width=0.75\linewidth]{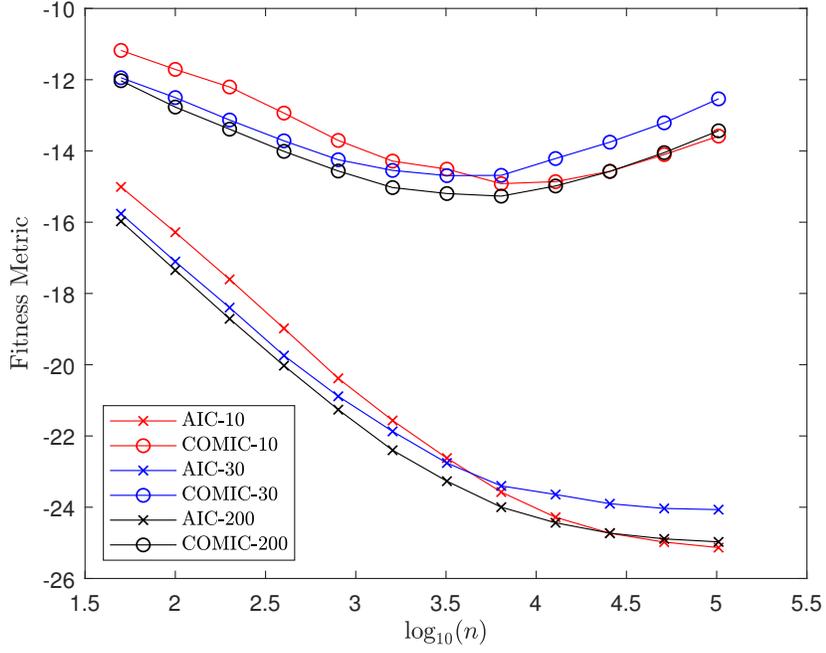}
    \caption{\rev{RWPT}}
    \label{fig:FMAICMultiPointsA}
\end{subfigure}
  \hspace*{-0.25in}
\begin{subfigure}[b]{\textwidth}
    \centering
    \includegraphics[width=0.75\linewidth]{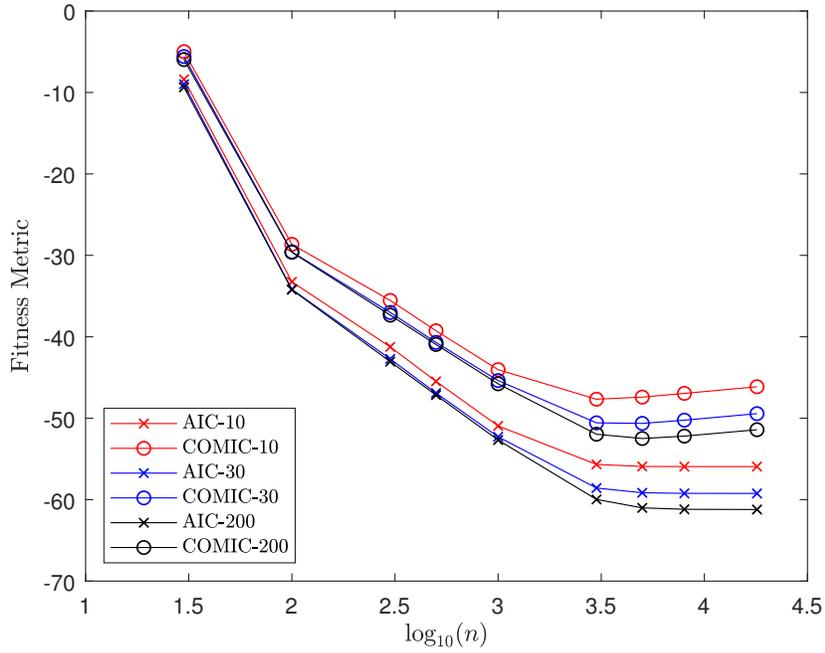}
    \caption{\rev{MTPT}}
    \label{fig:FMAICMultiPointsB}
\end{subfigure}
\caption{Fitness metrics for uniformly spaced data using $k=$ 10, 30, and 200 data points: RWPT (top) and MTPT (bottom).\\
}
\label{fig:FMAICMultiPoints}
\end{figure}

One problem that may arise when collecting data is a limited number of accessible locations.  The sparsity of data could change the number of particles needed to fit a given data set (and subsequently predict others).  We presume that this will also affect the optimal number of particles predicted by the COMIC. Thus, instead of assuming a large collection of data (e.g., $\geq 30$ points), we will also consider relatively small data sets ($\leq 10$ points) within the domain, to assess the degree to which parameter estimates are affected by less data and to what extent the optimal number of particles changes in order to achieve more precise estimates. Because we are relying on smaller collections of data, modification to the AIC is necessary \rev{\cite{CAVANAUGH1997}}. This corrected fitness metric is often denoted as the AICc and defined by \rev{\cite{HT1989}}
\begin{equation}
\text{AICc} = \text{AIC} + \dfrac{2p^2 + 2p}{k - p -1},
\end{equation}
where $k$ is the number of data points and $p$ is the number of parameters in the model. This leads to the obvious modification of the COMIC with $\Delta V = |\Omega|/n$ as
\begin{equation}
\text{COMICc} = \text{AICc} + \ln(n).
\end{equation}
Such a correction is relevant when we are attempting to compare different models with various numbers of data points and parameters.
\rev{Simulations were conducted by setting the final time $T = 1$ with $D=1$, $v = 0$ (i.e., no adjustable model parameters) and selecting $k=10$, $30$, or $200$ ``sample'' data points $\hat{c}_1, ..., \hat{c}_{k}$ from exact values of the analytic solution of the ADE \eqref{eq:ADE} on the interval $\Omega = [-5,5]$ with uniform spacing. We calculate the AICc and COMICc \rev{with different particle numbers} for both the random walk and mass transfer methods.} For the former, the optimal number of particles as given by the COMICc appears to increase for the sparse ($k=10$) data case to more than $n = 10^{4.3} \approx$ 20,000 particles (Fig. \ref{fig:FMAICMultiPointsA}).  On the other hand, the MTPT simulations appear relatively stable in terms of their COMIC fitness with an optimal number of particles \rev{around} $n =$ 3,000 (Fig. \ref{fig:FMAICMultiPointsB}).

\begin{figure}
\centering
\includegraphics[width=0.95\textwidth]{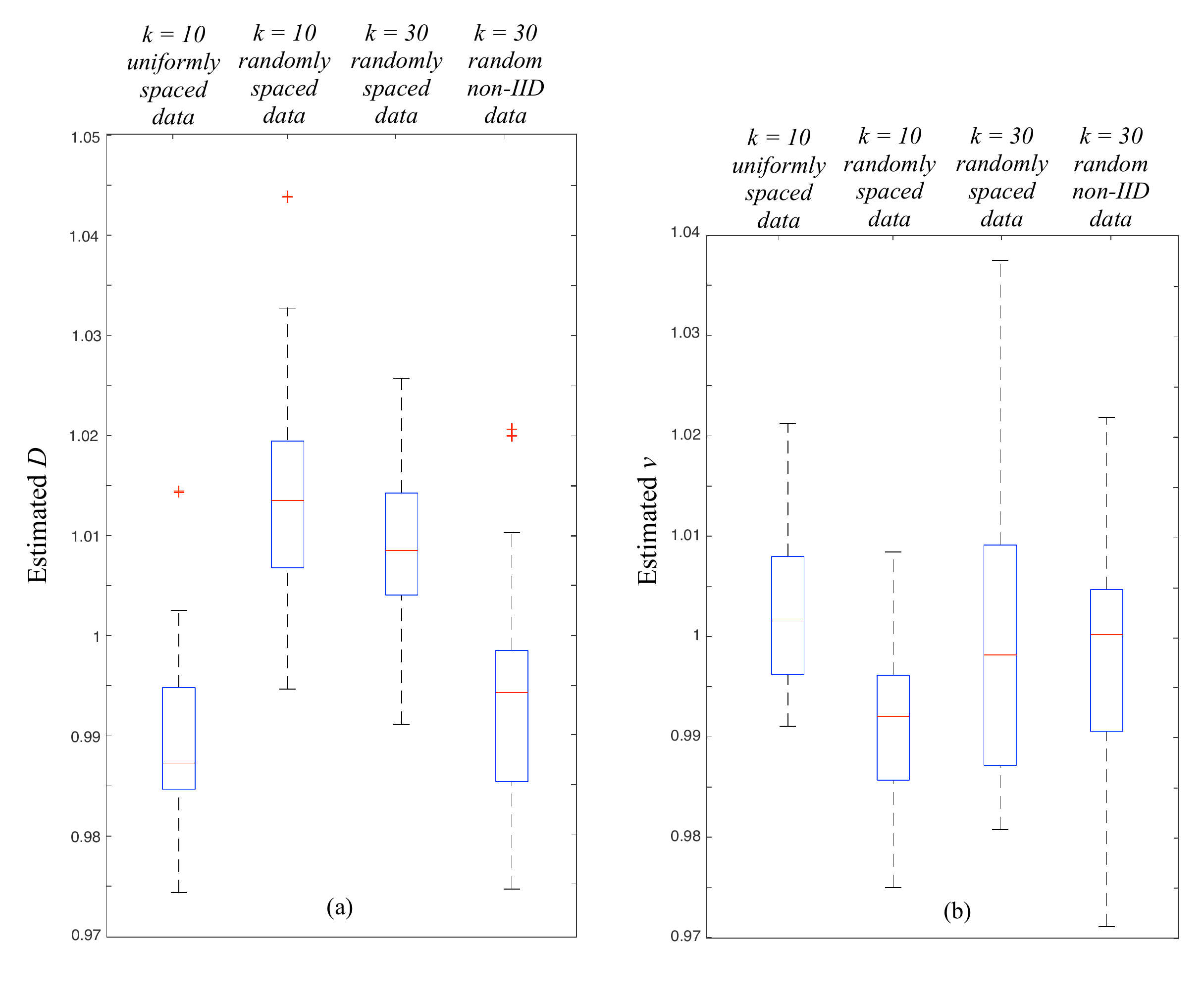}
\caption{Estimates of $D$ and $v$ using RWPT with different data spacings schemes, numbers, and distributional qualities.}
\label{fig:RWDVEst_all}
\end{figure}

We now use these results to fix the COMIC-optimal number of particles in simulations (20,000 for RWPT and 3,000 for MTPT) that now seek to estimate the parameters within the ADE. We choose data points from the analytic solution using $v=D=1$ and use MATLAB's built-in \texttt{fminsearch} function to minimize the AIC with initial guesses for both coefficients of $0.5$. Because the particles do not move in the MTPT algorithm, the solutions are the same if run multiple times (i.e., non-stochastic results). For 30 data points, the estimated values of \rev{$D$ and $v$ in the MTPT are within $10^{-6}$ of their true value $1$}.   On the other hand, the RWPT method gives a different result for each realization because of the random walks, so we show a box/whisker plot for those results in Fig. \ref{fig:RWDVEst_all}.

\begin{figure}
\centering
\includegraphics[width=0.75\textwidth]{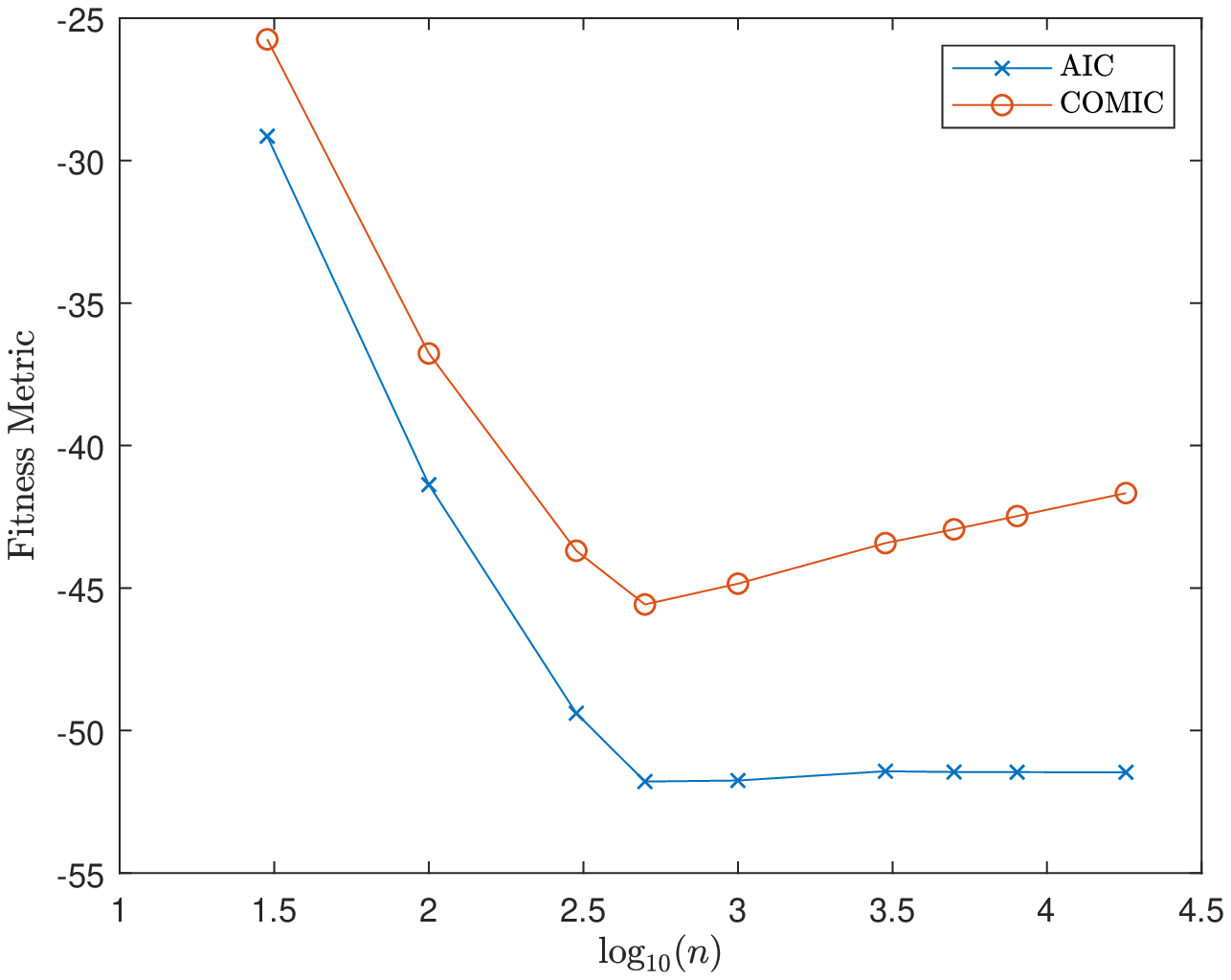}
\caption{Mass transfer fitness metrics \rev{determined by estimating the parameters $D$ and $v$ for} various particle numbers.}
\label{fig:MTAICDVEst}
\end{figure}

For even fewer data points (say, $k=6$), the random walk method requires significantly more particles to achieve a sufficient balance between goodness of fit and computational complexity. The qualitative behavior and shape of the AIC and COMIC for differing particle numbers are similar to Fig.~\ref{fig:FMAICMultiPoints}, and hence are not shown.
The mass transfer method displays similar behavior to other simulations - the COMIC optimal number of particles occurs at $n_{\text{COMIC}} \approx$ 3,000 and the estimated values of $D$ and $v$ are both \rev{within $10^{-6}$ of their true value of $1$}.
Hence, the MTPT is significantly more robust with respect to the size of the data.
Finally, we briefly mention that simulations were performed for every chosen particle number in which parameter values were first estimated and then used to compute the fitness metric. Again, the \rev{convex shape} of the AIC and COMIC for the RWPT method were nearly identical to Fig.~\ref{fig:FMAICMultiPoints}, though the new optimal number of particles implied by the COMIC for the MTPT is $n_{\text{COMIC}} \approx 500$. One reason for this decrease in optimal particle number is that the MTPT estimates are very close to the exact solution, so that any noise in the parameter estimation dominates the overall error. Figure \ref{fig:MTAICDVEst} displays the fitness metrics from MTPT simulations with $30$ data points. Using this new optimal particle number, we estimated the values of $D$ and $v$ to be \rev{within $10^{-5}$ of their true values, and these} are similar to the parameter estimates obtained using $3,000$ particles.

\section{Non-uniformly Spaced Data}
\label{sec:nonuniform}
Another issue that may arise under less ideal data collection circumstances is that observations from a field site may not be constrained to a certain grid, meaning that the data may not exist at uniformly-spaced gridpoints. Therefore, we examine the effects of such spatial heterogeneity on the COMIC and elucidate how this will affect parameter estimation. To generate non-uniform data, we randomly select data points within the domain of interest. More specifically, we perform simulations for $k=10$ and $k=30$ data points, with $D=1$, $v = 1$, $x_0 =0$ and $T = 1$. Additionally, the numerical spacing $\Delta V = 10/n$ remains constant. Calculations of the COMIC for random data display similar results to the case of equally-spaced data, so we will use the optimal number of particles from the previous section, as well. For the RWPT method, we perform simulations with 5,000 particles for $30$ data points and 20,000 particles for $10$ data points. For the MTPT method both simulations use 3,000 particles. The initial guess for both parameters is $0.5$ within all simulations. Due to the added randomness in the spacing of the data, we might anticipate that the parameter estimates would vary more than in the case of uniformly-spaced data; however, the results for these random walk parameter-estimation simulations are similar in both the magnitude and the variability of the estimates of both $D$ and $v$ (Fig. \ref{fig:RWDVEst_all}). For the MTPT method, a simulation with $30$ data points yields \rev{estimates of $D$ and $v$ within $10^{-5}$ of their true values}. Similarly, a simulation with $10$ data points provides \rev{estimates of $D$ and $v$ within $10^{-4}$ of their true values}. From these simulations, we conclude that the COMIC provides a useful and informative guide for the choice of particle number and parameter estimation even when the data is not uniformly spaced.
\rev{Next, we will consider alterations to the COMIC in order to address non-constant sampling volumes.}

\section{Non-Uniform Numerical Discretization Volume}
\label{sec:dv}
In the original calculation of the COMIC, the spatial volume $\Delta V(x)$ along which the solution is calculated is assumed to be a constant, and this results in equation \eqref{eq:COMIC} representing the computational fitness metric. However, in performing certain simulations, this assumption may not hold as either the binning of RWPT particles --- or positions of MTPT particles --- may not be evenly distributed throughout the physical domain. In contrast to the previous section, this issue only arises due to the computational method rather than the collected data. In such a case, $\Delta V(x)$ will vary with $x$, and the calculation of the COMIC must instead use equation \eqref{eq:COMIC_exact}. One immediate drawback from this formulation is that the value of the true concentration $c(x,t)$ at any point is unknown. Still, the numerical approximation $c_n(x,t)$ can be used to estimate this function within the COMIC, which becomes
\begin{equation}
\label{eq:COMIC_approx}
    \mathrm{COMIC} = \mathrm{AIC}  - \int c_n(x,t)\ln(\Delta V(x))\; dx.
\end{equation}


In this direction, we perform MTPT simulations in which the particles are initially placed randomly within the domain, and the random spacing between particles will determine the sampling volume $\Delta V(x)$ at each of the $n$ particle locations.
\begin{figure}
    \centering
    \includegraphics[width =0.85\textwidth]{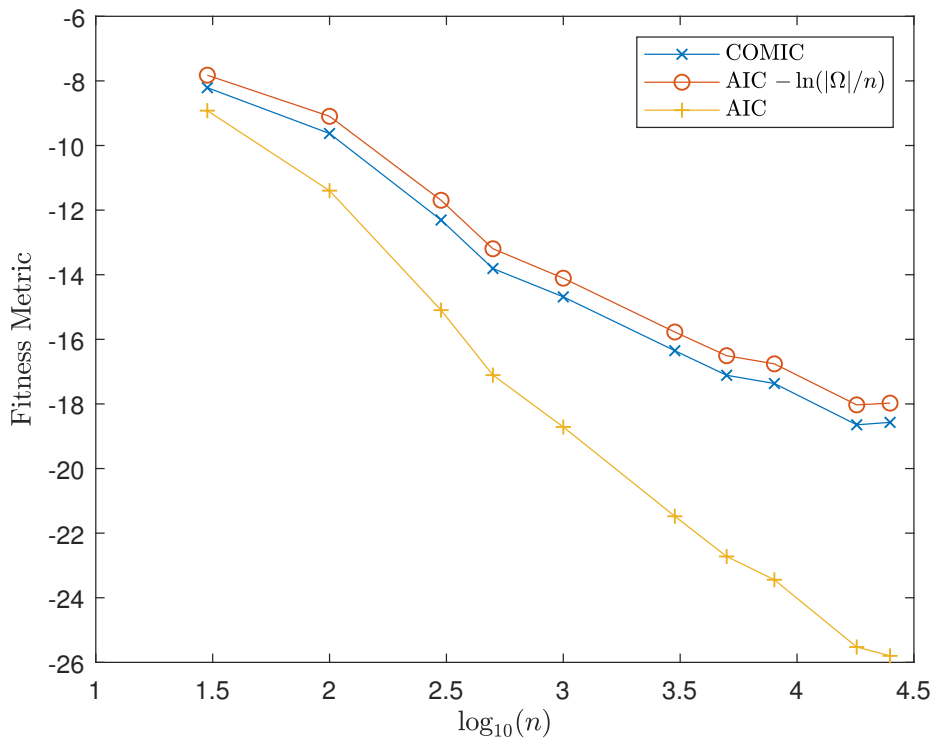}
    \caption{Fitness metrics versus particle number in the MTPT method with spatially-varying sampling volume $\Delta V(x)$}
    \label{fig:NonBin}
\end{figure}
This is performed with varying particle numbers, and the resulting COMIC of equation \eqref{eq:COMIC_approx} is computed for each simulation, once again using $k=30$ data points generated from the analytic solution of the ADE. To account for the randomness in the initial particle spacing, we compute the ensemble average of $30$ realizations (Fig. \ref{fig:NonBin}).
From these simulations, two key observations become apparent:
\begin{enumerate}
    \item The optimal number of MTPT particles given by the COMIC is $18,000$, which is much larger than the $3,000$ particles predicted by simulations with constant $\Delta V$ - see Fig.~\ref{fig:FMAICMultiPointsB}. This can be explained by the the randomness of the particle spacing. For instance, {\em Schmidt et al.} \cite{Schmidt_accuracy} show that the MTPT algorithm incurs increased error for randomly spaced, immobile particles because of mass-transfer ``gaps'' in areas of sparsely-distributed particles. The fitness of the method drops sharply compared to the uniformly-spaced scenario (Fig. \ref{fig:FMAICMultiPoints}), especially when the number of particles becomes small. This is clearly demonstrated by the AIC in Figure \ref{fig:NonBin}, as this curve possesses a large negative slope between $100$ and $18,000$ particles and only plateaus after $18,000$ particles. At this point the particles are distributed tightly enough in the domain to maintain the diffusion process regardless of the variability in their spacing. It should be noted that partitioning a portion of the diffusion process to random walks alleviates the problem of persistent mass-transfer ``gaps'' \cite{Entropy}.
    \item The exact COMIC calculation, which is performed by approximating the integral in \eqref{eq:COMIC_approx}, and the uniformly-spaced COMIC, given by equation \eqref{eq:COMIC}, have the same \rev{convex shape}.  For stationary particles, integrating a nonuniform $\Delta V(x)$ simply adds a constant value to computing the COMIC for constant $\Delta V$. Hence, the curve is simply shifted vertically, without further influencing its shape or the position at which it attains its minimum value, and this property persists throughout multiple simulations.  Both fitness metrics provide the identical optimal number of particles of $18,000$.  Thus, any calculation of the COMIC could use the equation \eqref{eq:COMIC} regardless of the sampling volume.
\end{enumerate}
\rev{In addition to changes in the COMIC that arise from data collection or choice of numerical method, we will also consider differing statistical assumptions that affect the information criterion, and this is performed in the next section.}

\section{Non-IID/Non-Gaussian Error Processes}
\label{sec:noGauss}
The original study of the COMIC \cite{Entropy} assumed that the residual between the (unknown) true concentration $c(x,t)$ and the data $\hat{c}_i$ for $i = 1, ..., k$ at each location was an independent and identically distributed (IID) Normal random variable, so that
\begin{equation}
c(x_i, T) = \hat{c}_i + \epsilon_i
\end{equation}
for every $i=1,...,k$ where $\epsilon_i \sim \mathcal{N}(0, \sigma^2)$ and the variance $\sigma^2$ must be estimated from the data.  This led to a log-likelihood function (hence AIC) that used the average sum of squared errors for the variance estimate within the fitness metric.  The assumption of IID Gaussian errors is not generally valid, and the properties of the error distribution are unknown in most cases. Here, we exploit the fact that the computational solution arises from a particle method in order to approximate the residuals.
In this scenario, \emph{Chakraborty et al} \cite{Chakraborty} showed that (a) the concentration approximation generated by any particle method is proportional to a binomial random variable that is only asymptotically Normal (as $n \to \infty$, $\Delta x \to 0$, and $\sqrt{n}\Delta x \to \infty$ where $\Delta x$ is the bin size of the method),  and (b) individual concentration errors could be treated as independent with variance proportional to the concentration, i.e. $\sigma^2_i = \dfrac{m_{\text{total}}}{n \Delta x} c_n(x_i,t)$, where $m_{\text{total}}$ is the total mass.
%
Indeed, this can be numerically verified using the RWPT and MTPT methods \cite{thesis}.
In \cite{Chakraborty}, an alternative information criterion was also proposed for selecting a ``best'' model over all parameter choices with the desirable property that the chosen parameter estimate $\hat{\theta}$ serves as a consistent estimator for the true parameter values $\theta$. More specifically, this criterion arises from an optimal fitting procedure that serves to minimize the weighted mean square error function
\begin{equation}
\label{mincrit}
\mcE(\theta) = \dfrac{1}{k} \sum_{i=1}^k w_i |\hat{c_i} - c_n(x_i,T; \theta)|^2
\end{equation}
where $\theta$ is the vector of unknown model parameters, the minimization weights are
\begin{equation}
w_i = \dfrac{1}{m_{\text{total}} \hat{c_i}},
\end{equation}
and the estimator $\hat{\theta}$ is given by
\begin{equation}
\hat{\theta} = \argmin_{\theta \in \mathbb{R}^p} \ \mcE(\theta).
\end{equation}
For the ADE problem described in previous sections, we merely have $\theta = [v, D]^T$.
We note that in the case that errors are normally-distributed, minimizing \eqref{mincrit} is equivalent to maximizing the log-likelihood function 
for a multivariate Gaussian distribution with $\hat{\sigma}_i^2 = m_{\text{total}} \hat{c_i}$ for $i = 1, \dots, N$ (see Appendix).

\begin{figure}
\centering
\hspace{-0.35in}
\begin{subfigure}[b]{\textwidth}
    \centering
    \includegraphics[width=0.75\textwidth]{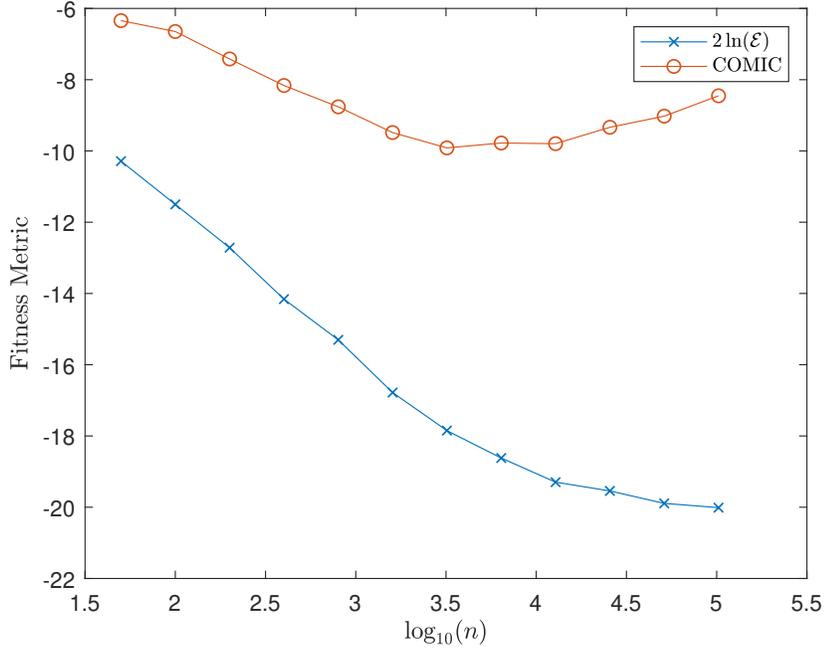}
    \caption{\rev{RWPT}}
    \label{fig:RWMTAICMGA}
\end{subfigure}
  \hspace*{-0.25in}
\begin{subfigure}[b]{\textwidth}
    \centering
    \includegraphics[width=0.75\linewidth]{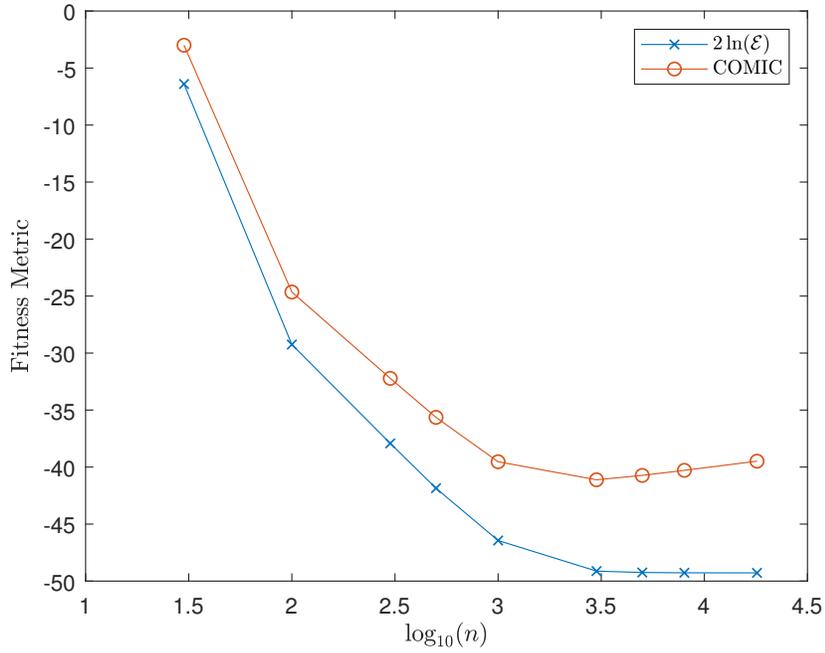}
    \caption{\rev{MTPT}}
    \label{fig:RWMTAICMGB}
\end{subfigure}
\caption{Fitness metrics for non-Gaussian error processes - RWPT (top) and MTPT (bottom). Here, the COMIC is computed from equation \eqref{COMICE}, and these curves can be directly compared with those of 
  Fig.~\ref{fig:FMAICMultiPoints} with $k=30$.\\
}
\label{fig:RWMTAICMG}
\end{figure}



As for the AIC, this information criterion does not account for the additional information incurred by taking large numbers of particles, and hence we augment it to create a new computational information criterion. \rev{Therefore}, in this case we the define the COMIC by
\begin{equation}
    \mathrm{COMIC} = 2\ln \left (\dfrac{1}{k} \sum_{i=1}^k \dfrac{1}{m_{\text{total}} \hat{c_i}} |\hat{c_i} - c_n(x_i,T)|^2 \right) + 2p - \int c_n(x,T)\ln(\Delta V(x))\; dx.
\end{equation}
where $\hat{c_i}$ is the concentration data at location $x_i$ and final time $t = T$. From the results of the previous section, it is beneficial to take $\Delta V$ constant, which reduces this formulation to
\begin{equation}
\label{COMICE}
    \mathrm{COMIC} = 2\ln \left (\dfrac{1}{k} \sum_{i=1}^k \dfrac{1}{m_{\text{total}} \hat{c_i}} (\hat{c_i} - c_n(x_i,T))^2 \right)+2p  - \ln(\Delta V).
\end{equation}
Using this criterion, we perform simulations of the random walk and mass transfer methods to compute the value of $2\ln(\mcE)$ and the COMIC. The formulation and implementation of these methods is analogous to that of the previous section with $k=30$ randomly-spaced data points.
From Figure \ref{fig:RWMTAICMG}, the optimal number of particles for the MTPT method is $3,000$, while for the RWPT method, it ranges between $3,200$ and $6,400$ particles. We take the midpoint of these and assume the optimal number of particles is about $5,000$. Then, this predicted value is essentially the same as that stemming from IID Gaussian error simulations \cite{Entropy}. Using the optimal number of particles to perform parameter estimation provides MTPT estimates \rev{within $10^{-4}$ of the true values of $D = 1$ and $v=1$}. Because the data is random, simulation outcomes may vary, but multiple runs with different data display similar qualitative behavior.  Estimated values of $D$ and $v$ using the RWPT method are also similar - see Figure \ref{fig:RWDVEst_all}. The maximal absolute error in the \rev{estimates of $D$ and $v$ are about $2.5\%$ and $3\%$, respectively, which are comparable to the parameter estimation performed for IID Gaussian error simulations}. Therefore, the COMIC demonstrates consistency among different error assumptions and estimators. In other words, optimal model discretization does not strongly depend on the weights in the sum of squared errors.


\begin{figure}
\centering
\hspace{-0.35in}
\begin{subfigure}[b]{\textwidth}
    \centering
  \includegraphics[width=0.75\textwidth]{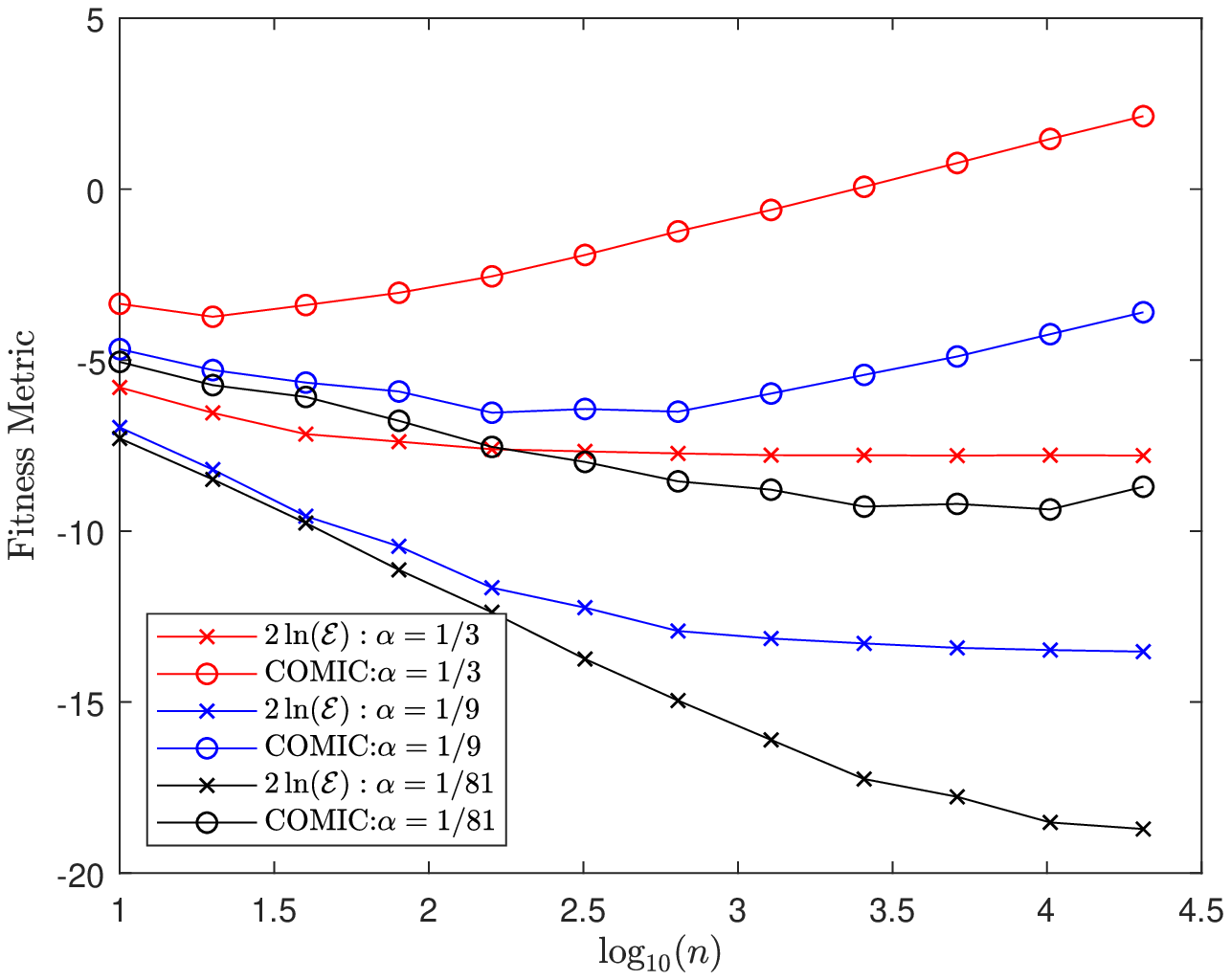}
    \caption{\rev{RWPT}}
    \label{fig:MTNonGauA}
\end{subfigure}
  \hspace*{-0.25in}
\begin{subfigure}[b]{\textwidth}
    \centering
      \includegraphics[width =0.75\textwidth]{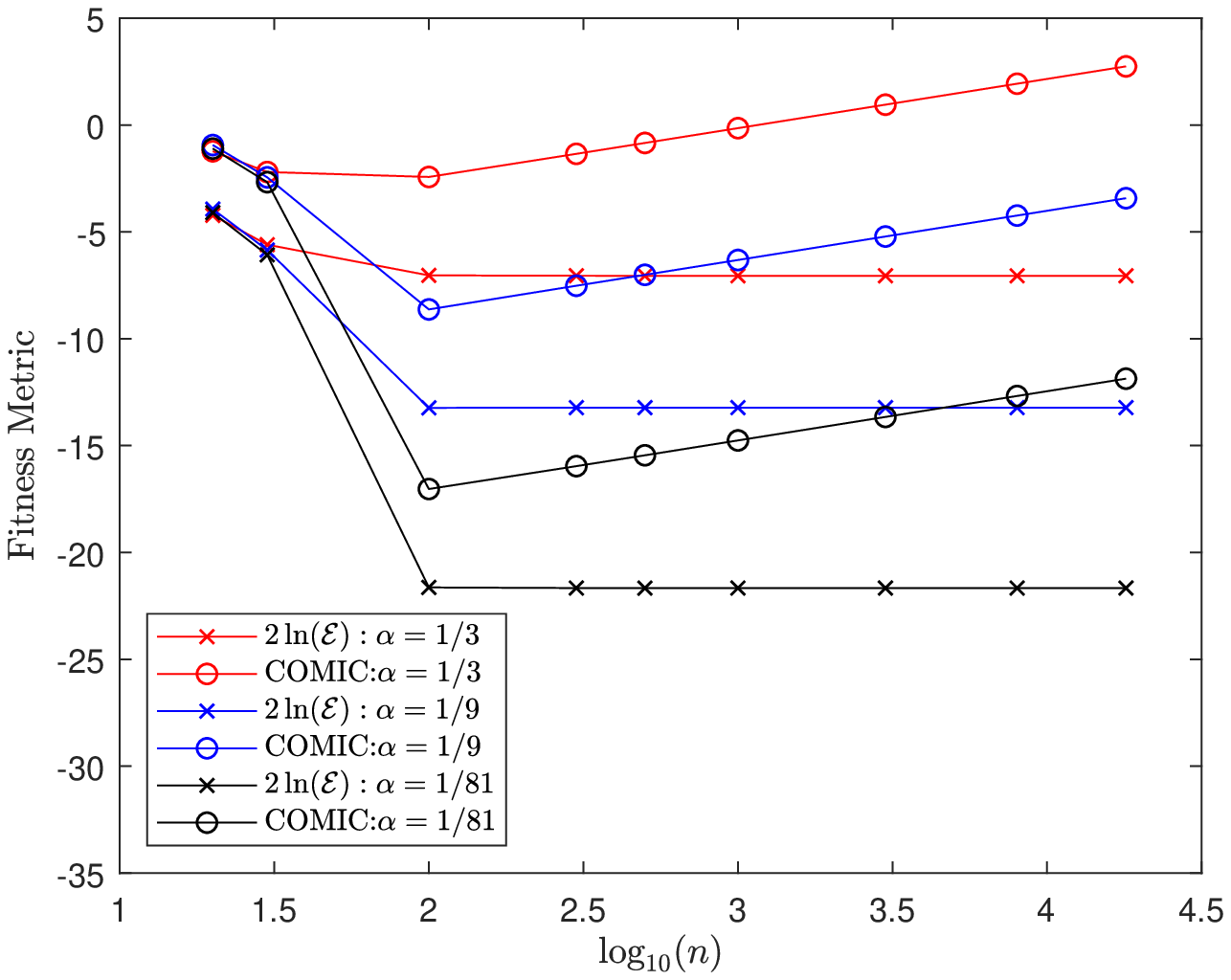}
    \caption{\rev{MTPT}}
    \label{fig:MTNonGauB}
\end{subfigure}
    \caption{Fitness metrics for Non-Gaussian Error Processes with noisy data - RWPT (top) and MTPT (bottom).\\
    }
    \label{fig:MTNonGau}
\end{figure}



Lastly, we perform similar simulations with noisy data, i.e. data that differs from the analytic solution due to measurement error modeled by independent random concentration noise at each spatial location.
Figure \ref{fig:MTNonGau} shows the fitness metrics arising from simulated data that is normally distributed with mean equal to the exact concentration from the analytic solution of eq. \eqref{eq:ADE} and standard deviation given by $\alpha = 1/3$, $1/9$, or $1/81$  of the exact concentration so that $\mathrm{Var} (\hat{c_i}) = \alpha^2 c(x_i,T)^2$. \rev{We note that any negative sample points are discarded, in order to guarantee non-negativity of the concentration data set.} Of course, the exact value of $\alpha$ can be tuned to adjust for the noise in the data.

There are several noteworthy features shown in a single realization of the simulations (different simulations show similar results):
\begin{enumerate}
   \item For the RWPT simulations, the optimal number of particles is strongly determined by the standard deviation of the independent noise held within the data.  In short, the noisier the data, the more noise can be contained within the model solution \rev{(i.e., fewer particles should be used)}.
    \item For the MTPT solution, the optimal value of the COMIC occurs at $n=100$ particles, which is far fewer than the optimal number of particles when $\alpha  = 0$.  In contrast to the RWPT solutions, the number of particles is not strongly affected by the magnitude of noise in the data.  This can be explained by the formulation of the COMIC, in which the dominant term in the calculation is the AIC, until $n$ becomes sufficiently large. The numerical simulation is converging to the exact solution of the PDE, i.e. a normal distribution, and as the number of particles nears $100$, the approximate solution closely agrees with the exact solution. However, assuming substantial noise in the data, the exact solution is far from the expected normal distribution, which means that the AIC will remain large even as the computational approximation converges. This can be seen by the graph of $2\ln(\mathcal{E})$ within Fig.~\ref{fig:MTNonGau}, which plateaus for simulations with more than $100$ particles.
    \item If a given data set is particularly noisy or fails to faithfully represent the underlying solution of the PDE, one would actually do best to use fewer particles in a simulation. When significant noise exists in the data, a low-noise, high-$n$ solution is not expected to be a better predictor because it is apt to be over-fit to peculiarities in a single noisy data set.
\end{enumerate}

\section{Conclusions}
\label{sec:conclusions}
%
We have investigated the use of the COMIC to select parsimonious and robust computational models for simple advective-diffusive transport in a variety of realistic, data-driven scenarios, including noisy, sparse, and spatially-heterogeneous data sets, non-uniform sampling volumes, and non-IID errors.
In the case of non-uniform sampling volumes, we have shown that the calculation of the COMIC can be further simplified using the average spacing between particles throughout the domain, and this does not influence the particle number selected by the algorithm.
The results of our simulations demonstrate that under any of these conditions, the COMIC is a flexible criterion that allows the user to select an appropriate number of particles in a simulation so as to guarantee the use of minimal computational information to construct a descriptive model.
In particular, we find that the use of large particles numbers is often superfluous in these simulations and needlessly increases the complexity of a model.
This highlights the importance of selecting a suitably efficient computational model with minimal information content to best make predictions based on a single given data set.
In particular, we find the following general rules regarding the optimal discretization for RWPT and MTPT simulations and/or parameter estimation of advection and diffusion in 1-D:
\begin{enumerate}

\item{For ``perfect'' (no noise) data that is uniformly spaced, the number of data does not strongly effect the optimal number of particles.  In general, RWPT (with binning) required $\approx 20,000$ particles, \rev{PTMT} required $\approx 3,000$.   
}
\item{For either solution technique, the parameters can be estimated quite accurately.  The optimal number of particles in the MTPT simulations drops to $\approx 500$ when estimating \rev{the velocity and diffusion parameters}.
}
\item{Non-uniformly spaced data does not change the ability to estimate parameters.
}
\item{Non-uniformly spaced particles in the MTPT technique significantly increases the optimal number of particles because of increased error in the simulations.
}
\item{Non-uniform numerical discretization adds a constant to the COMIC but does not change the shape of the COMIC versus discretization curve, hence the minimum is not shifted.  This means that the computational entropy penalty can be estimated by $\ln(n)$.
}
\item{The classical MLE fitness metric of average squared error, which comes from an assumption of IID Normal errors, is a reasonably good estimator for errors in which concentration variance is proportional to concentration.  However, as the noise in data increases, the optimal number of RWPT particles decreases sharply. In short, the error of the solution is directly tied to particle numbers, and hyper-accurate solutions are not representative of the noisy data (and may overfit the errors).  The accuracy of the MTPT method is not tied to particle numbers in the same way (once a minimum particle spacing is achieved), and so the optimal particle number remains $\approx 100$ for all noise levels.
}

\end{enumerate}

\newpage
\bibliography{biblio}

\appendix
\section{MLEs for non-IID Multivariate Gaussian Processes}

We first recall the process of obtaining maximum-likelihood parameter estimates $\hat{\theta}$ under the assumption that the errors between model and observations are independent, zero-mean Gaussians. In this case the likelihood function is given by
\bl
\be
L(y; \theta)=\left [(2\pi)^k |\Sigma(\theta)| \right]^{-1/2}\exp\left (-\frac{1}{2}y^T\Sigma(\theta)^{-1}y \right),
\ee
\el
where $k$ is the number of observation points, $\Sigma(\theta)$ is a diagonal (due to independence) covariance matrix of errors, and $y$ is a vector of residuals satisfying $y_i = \hat{c}_i-c(x_i, T)$ for $i=1,...,k$.
Now, if the errors are further assumed to be identically-distributed, then $\Sigma$ is a constant multiple of the identity and depends only upon a single variance parameter. In this case, $\Sigma = \sigma^2\mathbb{I}$ where $\sigma^2$ is the assumed variance of the error at each spatial data point $x_i$ for $i = 1, ..., k$.
The log-likelihood function then becomes
\bl
\be
\ln(L)=-\frac{k}{2}\ln(2\pi)-\frac{k}{2}\ln \sigma^2- \frac{k}{2\sigma^2} \frac{\SSE}{k}
\ee
\el
where $\SSE=y\cdot y = \vert y \vert^2$ represents the sum of squared errors.
Maximizing this function using standard tools provides an estimator of the observation variance, namely $\hat{\sigma}^2=\SSE/k$. Removing any constant terms that do not change from one model to another,  the corresponding log-likelihood evaluated at the MLE is
\bl
\be \label{SSEn}
\ln(\hat{L}) = - \ln \left ( \frac{\SSE}{k} \right ).
\ee
\el

Alternatively, if one does not assume that the errors are identically-distributed, then $\Sigma$ is not a constant multiple of the identity, but merely diagonal, so that
$$\Sigma = \begin{bmatrix}
\sigma_1^2 & 0 & \cdots & 0\\
0 & \sigma_2^2 & \cdots & 0\\
0 & 0 & \ddots & 0\\
0 & 0 & \cdots & \sigma_k^2
\end{bmatrix}.$$
Hence, the error distribution at each spatial data point $x_i$ possesses a different variance $\sigma_i^2$.
In this case, the MLE can be generated using similar methods as above, and the resulting values of $\hat{\sigma}_i^2$ are given by the multivariate sum of squared errors
$$
\hat{\sigma}_i^2 = \frac{1}{N} \sum_{\ell = 1}^N | (y_\ell)_i |^2
$$
where $y_1, ..., y_N$ is a collection of residual vectors with each $y_\ell$ representing a single sample from all spatial points, namely
$$(y_\ell)_i = \hat{c}_i-c(x_i, T)$$
for the $\ell$th sample with $i = 1, ..., k$.
Of course, if only one sample is collected at each spatial gridpoint so that $N=1$ and $y$ represents the single data vector, then the variance estimator reduces to $\hat{\sigma}_i^2 = y_i^2$.
Unfortunately, this does not provide a sharp estimate, as the standard deviation of the error distribution is merely equal to the data value at every spatial point.
Therefore, to provide a realistic estimate of these values, we would require multiple concentration measurements at each spatial data point and at a fixed time $T$. In the absence of this data or a suitable approximation, the variance of the error distribution cannot be accurately determined.
Because of this difficulty, \emph{Chakraborty et al} \cite{Chakraborty} proposed an alternative minimization criterion, which uses a weighted mean square error as described within Section \ref{sec:noGauss}, and this is derived only under the assumption that the error is approximately Gaussian.
This criterion enables us to generate a consistent parameter estimate when only a single concentration sample is available at each spatial location.

\end{document}